\documentclass[11pt]{amsart}
\usepackage[T1]{fontenc}
\usepackage{lmodern}
\usepackage{amssymb,mathtools}
\usepackage[letterpaper,hmargin=1.15in,top=1.1in,bottom=1.1in]{geometry}
\usepackage[hidelinks]{hyperref}

\newtheorem{theoremletter}{Theorem}

\newtheorem{lemma}{Lemma}[section]
\newtheorem{proposition}[lemma]{Proposition}
\newtheorem{corollary}[lemma]{Corollary}
\theoremstyle{remark}
\newtheorem{remark}[lemma]{Remark}
\numberwithin{equation}{section}

\newcommand{\N}{\mathbb{N}}
\newcommand{\Z}{\mathbb{Z}}
\newcommand{\Q}{\mathbb{Q}}
\newcommand{\T}{\mathsf{T}}
\newcommand{\crd}{\operatorname{crd}}
\newcommand{\rk}{\operatorname{rk}}
\newcommand{\detR}{\det\nolimits_R}
\newcommand{\detV}{\det\nolimits_V}
\newcommand{\symp}{\mathrm{sp}}
\newcommand{\Sp}{\mathrm{Sp}}

\begin{document}

\title[Cigler's quadratic recurrence for Hankel determinants]{A proof of Cigler's quadratic recurrence for Hankel determinants of Motzkin path generating functions}
\author{Tianxing Yu}
\address{Peking University, Beijing, China}
\email{2500015471@stu.pku.edu.cn}
\date{November 20, 2026}

\begin{abstract}
For an integer $k\ge0$ let $a_{n,k}(x,c)$ be the generating polynomial of the Motzkin paths from $(0,0)$ to $(n,k)$ that never go below height $0$, where level steps have weight $c+x$ at height $0$ and weight $c$ at positive height, and let $D_{k,n}(x,c)=\det\bigl(a_{i+j,k}(x,c)\bigr)_{0\le i,j\le n-1}$. We prove that for all $k\ge0$ and $n\ge2k+2$
\[
D_{k,n}D_{k,n-2k-2}-x^2D_{k,n-1}D_{k,n-2k-1}=(1-x^2)\,D_{k,n-k-1}^2,
\]
a quadratic recurrence conjectured by J.~Cigler in 2023 (eq.~(18) of arXiv:2309.15557 and MathOverflow question 456683). As a by-product we obtain an expression of $D_{k,n}$ as a $k\times k$ determinant built from Chebyshev polynomials reduced modulo $U_k-xU_{k-1}$, together with a bilinear identity for these $k\times k$ determinants that holds for all integer shifts. We also give a proof of Cigler's companion conjecture that $D_{k,n}(x,c)$ equals the determinant of an explicit matrix $H_{k,n}(x)$ with entries $0$, $1$ and powers of $x$, and in particular does not depend on $c$; this statement has recently been proved by V.~Jandhyala, and our proof was found independently. The proofs use Chebyshev polynomials, a model of $\Q(x)[X]/(U_k-xU_{k-1})$ by antisymmetric Laurent polynomials, a lemma on complementary minors, and a three-term Plücker relation.
\end{abstract}

\maketitle

\section{Introduction}

\subsection{Definitions and main results}
Let $x,c$ be indeterminates and let $k\ge0$ be an integer. Put $s_0=c+x$ and $s_j=c$ for $j\ge1$. Define $a_{n,j}=a_{n,j}(x,c)\in\Z[x,c]$ for $n\ge0$ and $j\in\Z$ by
\[
a_{0,j}=\delta_{j,0},\qquad a_{n,j}=0\quad(j<0),\qquad a_{n,j}=a_{n-1,j-1}+s_j\,a_{n-1,j}+a_{n-1,j+1}\quad(n\ge1,\ j\ge0).
\]
Thus $a_{n,j}$ is the weighted number of Motzkin paths from $(0,0)$ to $(n,j)$ that never go below height $0$, where up and down steps have weight $1$, level steps at height $0$ have weight $c+x$, and level steps at positive height have weight $c$. For $n\ge0$ let
\[
D_{k,n}=D_{k,n}(x,c)=\det\bigl(a_{i+j,k}\bigr)_{0\le i,j\le n-1},\qquad D_{k,0}=1.
\]
The main result of this paper is Theorem~\ref{thm:B} below, a quadratic recurrence for the sequence $(D_{k,n})_{n\ge0}$ conjectured by Cigler \cite{C23,CMO}.

For $p,q\ge0$ define
\[
h_{p,q}=h^{(k)}_{p,q}(x)=\begin{cases}
1 & \text{if } |p-q|\le k\le p+q \text{ and } p+q\equiv k\pmod 2,\\
x^{\,q-p-k} & \text{if } q-p>k,\\
0 & \text{otherwise},
\end{cases}
\]
and let $H_{k,n}(x)=\bigl(h_{p,q}\bigr)_{0\le p,q\le n-1}$ (the empty matrix if $n=0$) and $d_{k,n}(x)=\det H_{k,n}(x)$.
This is the matrix of Cigler \cite{CMO}, where its non-zero entries are described by $h_{k-m+l,\,m+l,\,k}=1$ for $0\le m\le k$, $l\ge0$, and $h_{l,\,k+i+l,\,k}=x^i$ for $i\ge1$, $l\ge0$: in the first clause $p=k-m+l$ and $q=m+l$, so that $p+q=k+2l\ge k$ and $q-p=2m-k\in[-k,k]$, in the second $q-p=k+i>k$. The condition $l\ge0$ is not written in \cite{CMO}, but it is forced by the example $H_{2,6}(x)$ given there (for $k=2$, $m=1$, $l=-1$ the first clause would give $h_{0,0}=1$).

\begin{theoremletter}[Cigler's first conjecture \cite{CMO}]\label{thm:A}
For all $k\ge0$ and $n\ge0$,
\[
D_{k,n}(x,c)=\det H_{k,n}(x)\quad\text{in }\Z[x,c].
\]
In particular, $D_{k,n}(x,c)$ does not depend on $c$.
\end{theoremletter}

A proof of Theorem~\ref{thm:A} has been given by Jandhyala \cite{J}. The proof in Appendix~\ref{app:A} was obtained independently; the two arguments are compared in \S\ref{ss:background}.

\begin{theoremletter}[Main theorem; Cigler's second conjecture \cite{CMO}, {\cite[eq.~(18)]{C23}}]\label{thm:B}
For all $k\ge0$ and all $n\ge2k+2$,
\begin{equation}\label{eq:main}
D_{k,n}\,D_{k,n-2k-2}-x^2\,D_{k,n-1}\,D_{k,n-2k-1}=(1-x^2)\bigl(D_{k,n-k-1}\bigr)^2\quad\text{in }\Z[x,c].
\end{equation}
\end{theoremletter}

By Theorem~\ref{thm:A}, the same identity holds for $d_{k,n}(x)=\det H_{k,n}(x)$; this is the form in which the recurrence is stated in \cite{CMO}.

In \cite{C23} the determinant $D_{k,n}(x,c)$ is denoted $D_{0,k,n}(b,c)$, with $b$ in place of $x$, and \eqref{eq:main} appears there as eq.~(18), without a range for $n$; the range $n\ge2k+2$ is the one given in \cite{CMO}.

A by-product of the proof (Lemma~\ref{lem:Lambda} and Proposition~\ref{prop:bilinear}) is a $k\times k$ determinant formula for $D_{k,n}$:

\begin{theoremletter}\label{thm:C}
Let $k\ge1$, let $U_m\in\Z[X]$ $(m\in\Z)$ be defined by $U_{-1}=0$, $U_0=1$ and $U_{m+1}=XU_m-U_{m-1}$ for all $m\in\Z$ $($so that $U_{-m-2}=-U_m)$, and let $P_k=U_k-xU_{k-1}$. The residue classes of $U_0,\dots,U_{k-1}$ form a basis of $\Q(x)[X]/(P_k)$. For $m\in\Z$ let $w(m)$ be the determinant of the $k\times k$ matrix whose rows are the coordinate vectors, with respect to this basis, of the residue classes of $U_{m-1},U_m,\dots,U_{m+k-2}$ modulo $P_k$. Since $P_k$ is monic in $X$ with coefficients in $\Z[x]$ $($and $U_j$ is monic in $X$ of degree $j)$, these coordinates lie in $\Z[x]$; hence $w(m)\in\Z[x]$. Then
\begin{gather}
D_{k,n}=(-1)^{nk+k(k+1)/2}\,w(n-k)\qquad\text{for all } n\ge0,\label{eq:C1}\\
w(\ell+k+1)\,w(\ell-k-1)-x^2\,w(\ell+k)\,w(\ell-k)=(1-x^2)\,w(\ell)^2\qquad\text{for all }\ell\in\Z.\label{eq:C2}
\end{gather}
\end{theoremletter}

The first formula \eqref{eq:C1} is, up to the passage from the residues of $P_n,\dots,P_{n+k-1}$ to those of $U_{n-k-1},\dots,U_{n-2}$, a special case of a known formula for Hankel determinants of linear combinations of moments \cite[Theorem~1]{K23}; see \S\ref{ss:background}. The second formula \eqref{eq:C2} is the new ingredient behind Theorem~\ref{thm:B}.

Theorems~\ref{thm:A} and~\ref{thm:B} are identities between polynomials in $x,c$ with integer coefficients, so they remain true after any specialization of $x$ and $c$ in any commutative ring.

\subsection{Background and previous work}\label{ss:background}
J.~Cigler posed Theorem~\ref{thm:A} as a conjecture in the MathOverflow question \cite{CMO} on 18 October 2023 and added the recurrence of Theorem~\ref{thm:B} to the question on 22 October 2023; the question also displays closed formulas for $k=1$ and $k=2$ (see Remark~\ref{rem:smallk}). They appear also in \cite{C23}, in the notation recalled after Theorem~\ref{thm:B}: there Cigler remarks that ``it seems that $D_{0,k,n}(b,c)$ does not depend on $c$ either'', and, after a list of linear recurrences in $n$ (eq.~(17) of \cite{C23}), that these sequences ``also seem to satisfy'' the quadratic recurrence of Theorem~\ref{thm:B} (eq.~(18) of \cite{C23}). The independence of $c$ had been conjectured earlier as eq.~(2.4) of \cite[Conjecture~2.2]{C22},\footnote{Conjecture~2.2 of \cite{C22} consists of two statements, and only the second one, eq.~(2.4), concerns the present paper. The first statement, as printed, fails already for $k=m=1$; it is not discussed here.} in the form $d^{(k)}_0(n,t+s,s)=d^{(k)}_0(n,t,0)$, which in our notation reads $D_{k,n}(-t,t+s)=D_{k,n}(-t,t)$ (see Remark~\ref{rem:smallk} for the translation).

\emph{Theorem~\ref{thm:A}.} The independence of $c$ alone is an instance of the classical invariance of Hankel determinants under the binomial transform. Write $J(c)=J(0)+cI$, where $J(c)$ is the matrix $J$ of Section~\ref{sec:prelim} and $J(0)$ is the same matrix for $c=0$. Then $a_{n,k}(x,c)=e_0^{\T}(J(0)+cI)^ne_k=\sum_j\binom nj c^{\,n-j}a_{j,k}(x,0)$, and therefore, by the Vandermonde convolution $\sum_{l+l'=m}\binom il\binom j{l'}=\binom{i+j}m$, the Hankel matrix $(a_{i+j,k}(x,c))_{0\le i,j\le n-1}$ is the Hankel matrix $(a_{i+j,k}(x,0))_{0\le i,j\le n-1}$ multiplied on the left by the lower unitriangular matrix $\bigl(\binom il c^{\,i-l}\bigr)_{0\le i,l\le n-1}$ and on the right by its transpose. See \cite{L} (for $c=1$; the computation is the same for general $c$), \cite[Theorem~6.1]{SS} (for general $c$) and \cite{P}; Cigler also asked about this invariance on Mathematics Stack Exchange \cite{CMSE}. In the present paper it is reproved, in a slightly different form, as Corollary~\ref{cor:cindep}.

The full statement of Theorem~\ref{thm:A} has been proved by V.~Jandhyala \cite{J}, in a preprint of 18 November 2026 that includes a formalization in Lean~4 over an arbitrary commutative ring. His argument removes $c$ by the binomial transform, cuts each path at an intermediate step, which turns the Hankel determinant into a determinant of counts of paths between two heights, and reduces the factorization of the latter matrix as $H_{k,n}(x)$ times a unitriangular matrix to a local identity for the entries of $H_{k,n}(x)$. The paper \cite{J} states explicitly that the quadratic recurrence remains open; it is checked there numerically for $k\le5$ and $n<22$. The proof of Theorem~\ref{thm:A} given in Appendix~\ref{app:A} was found independently. It also cuts the paths at an intermediate step, but it computes the weighted counts of paths between arbitrary heights (Lemma~\ref{lem:arb}) by means of the linearization formula for Chebyshev polynomials (Lemma~\ref{lem:cheb-lin}), an algebraic form of the reflection principle, and of the generating function of Lemma~\ref{lem:gf}, and it keeps $c$ throughout.

\emph{Theorem~\ref{thm:C}.} Since $D_{k,n}$ does not depend on $c$, we may take $c=0$. Then, by Lemma~\ref{lem:jacobi}(iii), $a_{n,k}=L(X^nP_k)$ is a linear combination, with leading coefficient $1$, of the $k+1$ consecutive moments $L(X^n),\dots,L(X^{n+k})$ of the functional $L$ of the weights $(x,0,0,\dots)$, whose monic orthogonal polynomials are the $P_m=U_m-xU_{m-1}$; thus $D_{k,n}$ is a Hankel determinant of a linear combination of moments. For such determinants \cite[Theorem~1]{K23} gives a determinant formula of size $d\times d$, where $d+1$ is the length of the linear combination. In our case $d=k$, the Hankel determinants of the moments $L(X^n)$ themselves are equal to $1$, and the formula reads
\[
D_{k,n}=(-1)^{nk}\,\frac{\det_{1\le i,j\le k}\bigl(P_{n+i-1}(\rho_j)\bigr)}{\prod_{1\le i<j\le k}(\rho_j-\rho_i)},
\]
where $\rho_1,\dots,\rho_k$ are the roots of $P_k$ in an algebraic closure of $\Q(x)$. When the $\rho_j$ are pairwise distinct, expressing $P_{n+i-1}(\rho_j)$ through the coordinates of the residue of $P_{n+i-1}$ modulo $P_k$ turns the right-hand side into the $k\times k$ determinant of Lemma~\ref{lem:reduction}; together with Lemma~\ref{lem:Lambda} this gives the first formula of Theorem~\ref{thm:C}. That formula is therefore a special case of \cite[Theorem~1]{K23}, a formula of Christoffel type which was stated by Elouafi \cite{E} and, in the language of random matrices, by Brézin and Hikami \cite{BH}; see \cite[Section~2]{K23} for its history. Our proof of Lemma~\ref{lem:reduction}, by complementary minors of the band matrix $P_k(J)$, is included to keep the paper self-contained; it needs no hypothesis on the roots of $P_k$. Closed formulas for Hankel determinants of linear combinations of three and of four consecutive moments, for weights that are constant above height zero (which is our situation), are given in \cite[Corollaries~3 and~6]{CK21}, and Chern and Shi \cite{CS} use \cite[Theorem~1]{K23} for $x\in\{0,1\}$. We have not found the second formula of Theorem~\ref{thm:C} in the literature.

\emph{Theorem~\ref{thm:B}.} Several special cases are known or follow directly from known results.
\begin{itemize}
\item For $k=0$ the statement is trivial, since $D_{0,n}=1$.
\item For $k=1$, $D_{1,n}$ is a Chebyshev (Fibonacci-type) polynomial (Remark~\ref{rem:smallk}), and Cigler notes in \cite{C23}, after eq.~(18), that in this case the recurrence follows from Cassini-type identities.
\item For $x=0$ the determinants are known in closed form: $D_{k,(k+1)m}(0,c)=(-1)^{\binom{k+1}2m}$ and $D_{k,n}(0,c)=0$ if $k+1\nmid n$ \cite[Theorem~1.1]{CS}, where this is derived from \cite{CK11}; the evaluation is also quoted in \cite{CMO}. The recurrence then reduces to $D_{k,n}D_{k,n-2k-2}=D_{k,n-k-1}^2$, which follows at once.
\item For $x=1$ the sequence $(D_{k,n}(1,c))_{n\ge0}$ is periodic with period $2k+1$ \cite[Theorem~1.6]{CS}, and the recurrence reduces to $D_{k,n}D_{k,n-2k-2}=D_{k,n-1}D_{k,n-2k-1}$, which again follows at once.
\end{itemize}
Beyond these cases we are not aware of results that imply Theorem~\ref{thm:B}. The paper \cite{CS} treats only the specializations $x\in\{0,1\}$ and contains neither the matrix $H_{k,n}(x)$ nor the recurrence for general $x$. For $k\ge1$, the follow-up \cite{C24} to \cite{C23} treats only constant sequences $s=(c,c,c,\dots)$. The paper \cite{CK21} gives the closed formulas mentioned above but not the quadratic recurrence. Finally, \cite[Corollary~9]{K23} gives a linear recurrence with constant coefficients of order $2^k$ for $n\mapsto D_{k,n}$, which is the statement (17) of \cite{C23} (see Remark~\ref{rem:smallk}); it does not yield the quadratic recurrence of Theorem~\ref{thm:B}. Theorem~\ref{thm:B} has the shape of a three-term bilinear recurrence of Gale--Robinson (Somos) type, with span $2k+2$ and shifts $1$ and $k+1$; recurrences of exactly this shape, with coefficients $1$ and $-1$, hold for the shifted Hankel determinants of the $q$-deformed metallic numbers. They were conjectured by Ovsienko and Pedon \cite[(1.22)]{OP} and proved by Han and Pedon \cite[Theorem~C]{HP}; those sequences are different from ours. To the best of our knowledge, as of 19 November 2026 neither a proof nor a counterexample of Theorem~\ref{thm:B} for general $x$ and $k\ge2$ has been made public, and the question \cite{CMO} has no answers.\footnote{Eq.~(18) of \cite{C23} was also listed as an open problem on 4 November 2026 in a public GitHub repository of automated proof attempts (\texttt{the-omega-institute/trureturing}, issue 26512), where it was closed as unresolved.}

\subsection{What is new}
To the best of our knowledge, the following results are new: Theorem~\ref{thm:B} for general $x$ and $k\ge2$, and the second formula of Theorem~\ref{thm:C}, that is, the bilinear identity of Proposition~\ref{prop:bilinear} for the $k\times k$ determinants $w(\ell)$, on which the proof of Theorem~\ref{thm:B} rests, together with the Laurent-polynomial model of \S\ref{ss:laurent} and the operator $\Lambda$ of Lemma~\ref{lem:Lambda} used in its proof. A shorter proof of Theorems~\ref{thm:B} and~\ref{thm:C}, in which \eqref{eq:C2} is obtained as a reduction of the Hirota--Miwa equation for Casorati determinants, is given in \cite{Sun}.

The other results of the paper are known or are consequences of known results, and are included to keep the paper self-contained. Theorem~\ref{thm:A} was first proved in \cite{J}; Appendix~\ref{app:A} gives an alternative proof, found independently. The first formula of Theorem~\ref{thm:C} is a special case of \cite[Theorem~1]{K23}, of which Lemma~\ref{lem:reduction} gives a root-free proof. The closed formulas for $k=1,2$ and the linear recurrences (17) of \cite{C23} follow from \cite[Corollary~3]{CK21} and \cite[Corollary~9]{K23}; Remark~\ref{rem:smallk} indicates how they also follow from Theorem~\ref{thm:C}. The ingredients of the proof of Theorem~\ref{thm:A} (the description $a_{n,j}=(J^n)_{0,j}$ by a transfer matrix, the linearization $U_pU_k=\sum_{i=0}^kU_{p-k+2i}$ of Chebyshev polynomials, and the generating function $z^q\Gamma^{q+1}$) are standard.

\enlargethispage{3pt}\emph{Idea of the proof of Theorem~\ref{thm:B}.} By Corollary~\ref{cor:cindep} we may take $c=0$; the proof of Theorem~\ref{thm:B} does not use Theorem~\ref{thm:A}, except for the reformulation in terms of $\det H_{k,n}(x)$. For $c=0$, Lemma~\ref{lem:jacobi} shows that $D_{k,n}$ is the $n$-th Hankel determinant of the functional $f\mapsto L(fP_k)$, a Christoffel-type modification of the moment functional $L$ of the weights $(x,0,0,\dots)$, whose monic orthogonal polynomials are $P_m=U_m-xU_{m-1}$. A lemma on complementary minors (Lemma~\ref{lem:compl}) reduces this $n\times n$ determinant to a $k\times k$ determinant formed from the residues of $P_n,\dots,P_{n+k-1}$ modulo $P_k$ (Lemma~\ref{lem:reduction}). The substitution $X=t+t^{-1}$ identifies $\Q(x)[X]/(P_k)$ with a space $V$ of antisymmetric Laurent polynomials modulo $Q(t)=t^{2k+2}-xt^{2k+1}+xt-1$. In $V$ the classes $F_j$ of $t^j-t^{-j}$ satisfy a four-term linear relation, and multiplication by a symmetric Laurent polynomial induces an operator $\Lambda$ on $V$ with $\det\Lambda=\pm1$ which shifts the index $j$ by $k+1$ in two different ways (Lemmas~\ref{lem:laurent} and~\ref{lem:Lambda}). This gives the first formula of Theorem~\ref{thm:C} and expresses the five determinants occurring in Theorem~\ref{thm:B} through the $k\times k$ minors of a single window of $k+2$ consecutive vectors $F_j$. If the window has rank $k$, its kernel is two-dimensional, and the three-term Plücker relation reduces the bilinear identity to an identity for two polynomials $\gamma(t),\eta(t)$ of degree at most $k+1$ spanning the kernel. The Laurent polynomial $\gamma(t)\eta(t^{-1})-\eta(t)\gamma(t^{-1})$ is divisible by $Q(t)$, and a comparison of exponents shows that it is a constant multiple of $t^{-k-1}Q(t)$. Evaluating this identity at $t=x$, where $Q(x)=x^2-1$, produces the factor $1-x^2$ in Theorem~\ref{thm:B}. Windows of rank less than $k$ do occur, also in the range used for Theorem~\ref{thm:B} (Remark~\ref{rem:degenerate}); for them all five determinants vanish.

\emph{Organization.} Section~\ref{sec:prelim} collects conventions and two elementary lemmas, Section~\ref{sec:proofBC} proves Theorems~\ref{thm:B} and~\ref{thm:C}, Section~\ref{sec:remarks} contains remarks (the cases $k=1,2$, degenerate windows, the operator $\Lambda$, a connection with symplectic characters, and computational checks), and Appendix~\ref{app:A} proves Theorem~\ref{thm:A}. Section~\ref{sec:proofBC} uses only Section~\ref{sec:prelim}; it uses Theorem~\ref{thm:A} only for the reformulation of Theorem~\ref{thm:B} for $d_{k,n}=\det H_{k,n}(x)$ stated after it.

\subsection{Notation}
$\N=\{0,1,2,\dots\}$; $\Z$ and $\Q$ are the integers and the rationals, and $A$ denotes a commutative ring. A binomial coefficient $\binom nj$ is $0$ unless $0\le j\le n$. For a power series $f$ in $z$, $[z^m]f$ is the coefficient of $z^m$ in $f$. We write $\delta_{i,j}$ for the Kronecker delta and $W^{\T}$ for the transpose of a matrix or vector $W$. All matrices and vectors are indexed from $0$.

\section{Preliminaries}\label{sec:prelim}

\subsection{Infinite matrices and the matrix $J$}
Let $A$ be a commutative ring. Let $M(A)$ be the set of matrices $W=(W_{i,j})_{i,j\in\N}$ with entries in $A$ such that every row and every column of $W$ has only finitely many non-zero entries. With the usual matrix product (all sums are finite), $M(A)$ is an associative ring with unit $I$; hence $f(W)\in M(A)$ is defined for $W\in M(A)$ and $f\in A[X]$, and polynomials in the same $W$ commute. Let $e_m$ $(m\in\N)$ be the $m$-th standard basis column vector (finitely supported); then $We_m$ is the $m$-th column of $W$ and $e_i^{\T}We_j=W_{i,j}$. If $W$ is symmetric, then so is $f(W)$, and $e_i^{\T}f(W)=\bigl(f(W)e_i\bigr)^{\T}$.

Let $J\in M(\Z[x,c])$ be the symmetric tridiagonal matrix with
\[
J_{m,m}=s_m\quad(s_0=c+x,\ s_m=c\text{ for }m\ge1),\qquad J_{m,m+1}=J_{m+1,m}=1,
\]
all other entries being $0$.

\begin{lemma}\label{lem:transfer}
For all $n,j\ge0$ we have $a_{n,j}=(J^n)_{0,j}$. Moreover $a_{n,j}=0$ for $j>n$ and $a_{n,n}=1$.
\end{lemma}

\begin{proof}
We have $(J^0)_{0,j}=\delta_{0,j}$, and for $n\ge1$
\[
(J^n)_{0,j}=\sum_i(J^{n-1})_{0,i}J_{i,j}=(J^{n-1})_{0,j-1}+s_j\,(J^{n-1})_{0,j}+(J^{n-1})_{0,j+1},
\]
where the first term is to be read as $0$ if $j=0$ (column $j$ of $J$ has the entries $J_{j-1,j}=1$ if $j\ge1$, $J_{j,j}=s_j$ and $J_{j+1,j}=1$). These are the recurrence and the initial condition defining $a_{n,j}$ $(j\ge0)$. The triangularity follows from the recurrence by induction on $n$: if $j>n$, then $j-1,j,j+1>n-1$, so $a_{n,j}=0$; and $a_{n,n}=a_{n-1,n-1}+0+0=1$.
\end{proof}

\subsection{Chebyshev polynomials}
Let $U_m\in\Z[X]$ $(m\in\Z)$ be defined, as in Theorem~\ref{thm:C}, by $U_{-1}=0$, $U_0=1$ and the two-sided recurrence $U_{m+1}=XU_m-U_{m-1}$ $(m\in\Z)$. Thus $U_1=X$ and $U_{-2}=XU_{-1}-U_0=-1$.

\begin{lemma}\label{lem:cheb}
In $\Z[t,t^{-1}]$ we have, for all $m\in\Z$,
\begin{equation}\label{eq:cheb}
(t-t^{-1})\,U_m(t+t^{-1})=t^{m+1}-t^{-m-1},
\end{equation}
and in $\Z[X]$
\begin{equation}\label{eq:negU}
U_{-m-2}=-U_m\qquad\text{for all }m\in\Z.
\end{equation}
\end{lemma}

\begin{proof}
Both sides of \eqref{eq:cheb}, viewed as sequences indexed by $m\in\Z$, satisfy the two-sided recurrence $f_{m+1}=(t+t^{-1})f_m-f_{m-1}$ (for the right-hand side, $(t+t^{-1})(t^{m+1}-t^{-m-1})-(t^m-t^{-m})=t^{m+2}-t^{-m-2}$), and they agree at $m=-1$ and $m=0$; hence they agree for all $m$. Similarly, the sequence $(-U_{-m-2})_{m\in\Z}$ satisfies the same two-sided recurrence as $(U_m)$ (multiply $U_{-m-1}=XU_{-m-2}-U_{-m-3}$ by $-1$) and agrees with it at $m=-1$ and $m=0$.
\end{proof}

\section{Proof of Theorems~\ref{thm:B} and~\ref{thm:C}}\label{sec:proofBC}

The case $k=0$ of Theorem~\ref{thm:B} is trivial. By Lemma~\ref{lem:transfer}, the numbers $a_{n,j}=a_{n,j}(x,c)$ are the numbers $a^{(c)}_{n,j}$ of Corollary~\ref{cor:cindep} below, so Lemma~\ref{lem:jacobi}(iv) below, applied with $s=(c+x,c,c,\dots)$ and $P^{(c)}_0=1$, gives $D_{0,n}=\det I_n=1$ (alternatively, $H_{0,n}$ is upper unitriangular and Theorem~\ref{thm:A} applies). The identity of Theorem~\ref{thm:B} then reads $1-x^2=1-x^2$. From now on $k\ge1$.

\subsection{$D_{k,n}$ as a principal minor}
We use the ring $M(A)$ of infinite matrices and the standard basis vectors $e_m$ introduced in Section~\ref{sec:prelim}.

\begin{lemma}\label{lem:jacobi}
Let $A$ be a commutative ring and let $s=(s_0,s_1,\dots)$ be an arbitrary sequence of elements of $A$. Let $J^{(s)}\in M(A)$ be the symmetric tridiagonal matrix with $J^{(s)}_{m,m}=s_m$, $J^{(s)}_{m,m+1}=J^{(s)}_{m+1,m}=1$ and all other entries $0$, and put $a^{(s)}_{n,j}:=\bigl((J^{(s)})^n\bigr)_{0,j}$ for $n,j\ge0$. $($By the argument in the proof of Lemma~\ref{lem:transfer}, which uses only the column structure of $J^{(s)}$, these elements satisfy $a^{(s)}_{0,j}=\delta_{j,0}$ and $a^{(s)}_{n,j}=a^{(s)}_{n-1,j-1}+s_j\,a^{(s)}_{n-1,j}+a^{(s)}_{n-1,j+1}$ for $n\ge1$, $j\ge0$, where $a^{(s)}_{n-1,-1}:=0$; this is the recurrence of Section~1 with the weights $s_j$.$)$ Define monic polynomials $P^{(s)}_m\in A[X]$ by
\[
P^{(s)}_0=1,\qquad P^{(s)}_1=X-s_0,\qquad P^{(s)}_{m+1}=(X-s_m)P^{(s)}_m-P^{(s)}_{m-1}\quad(m\ge1),
\]
and let $L^{(s)}\colon A[X]\to A$ be the $A$-linear map $L^{(s)}(f)=\bigl(f(J^{(s)})\bigr)_{0,0}=e_0^{\T}f(J^{(s)})e_0$. Then for all $i,j,k,m,n\ge0$:
\begin{enumerate}
\item[(i)] $P^{(s)}_m(J^{(s)})\,e_0=e_m$;
\item[(ii)] $\bigl(P^{(s)}_k(J^{(s)})\bigr)_{i,j}=L^{(s)}\bigl(P^{(s)}_iP^{(s)}_jP^{(s)}_k\bigr)$; in particular $(k=0)$ $L^{(s)}\bigl(P^{(s)}_iP^{(s)}_j\bigr)=\delta_{i,j}$;
\item[(iii)] $a^{(s)}_{n,k}=L^{(s)}\bigl(X^nP^{(s)}_k\bigr)$;
\item[(iv)] $\det\bigl(a^{(s)}_{i+j,k}\bigr)_{0\le i,j\le n-1}=\det\bigl((P^{(s)}_k(J^{(s)}))_{i,j}\bigr)_{0\le i,j\le n-1}$.
\end{enumerate}
\end{lemma}

\begin{proof}
(i) We have $J^{(s)}e_0=s_0e_0+e_1$ and $J^{(s)}e_m=e_{m-1}+s_me_m+e_{m+1}$ for $m\ge1$. So $P^{(s)}_1(J^{(s)})e_0=e_1$ and, inductively for $m\ge1$,
$P^{(s)}_{m+1}(J^{(s)})e_0=(J^{(s)}-s_mI)e_m-e_{m-1}=e_{m+1}$.

(ii) Using (i) and the symmetry of $P^{(s)}_i(J^{(s)})$,
\begin{align*}
\bigl(P^{(s)}_k(J^{(s)})\bigr)_{i,j}&=e_i^{\T}P^{(s)}_k(J^{(s)})e_j=\bigl(P^{(s)}_i(J^{(s)})e_0\bigr)^{\T}P^{(s)}_k(J^{(s)})P^{(s)}_j(J^{(s)})e_0\\
&=e_0^{\T}P^{(s)}_i(J^{(s)})P^{(s)}_k(J^{(s)})P^{(s)}_j(J^{(s)})e_0=L^{(s)}\bigl(P^{(s)}_iP^{(s)}_jP^{(s)}_k\bigr).
\end{align*}
For $k=0$, $\bigl(P^{(s)}_0(J^{(s)})\bigr)_{i,j}=I_{i,j}=\delta_{i,j}$.

(iii) $a^{(s)}_{n,k}=e_0^{\T}(J^{(s)})^ne_k=e_0^{\T}(J^{(s)})^nP^{(s)}_k(J^{(s)})e_0=L^{(s)}\bigl(X^nP^{(s)}_k\bigr)$.

(iv) Write $P^{(s)}_i=\sum_{l\le i}\xi_{i,l}X^l$ with $\xi_{i,i}=1$, and let $\Xi=(\xi_{i,l})_{0\le i,l\le n-1}$, a lower unitriangular matrix. By (ii), (iii) and the linearity of $L^{(s)}$,
\begin{align*}
\Bigl(\bigl(P^{(s)}_k(J^{(s)})\bigr)_{i,j}\Bigr)_{i,j<n}&=\Bigl(L^{(s)}\bigl(P^{(s)}_iP^{(s)}_jP^{(s)}_k\bigr)\Bigr)_{i,j<n}=\Xi\cdot\Bigl(L^{(s)}\bigl(X^{l+l'}P^{(s)}_k\bigr)\Bigr)_{l,l'<n}\cdot\Xi^{\T}\\
&=\Xi\cdot\bigl(a^{(s)}_{l+l',k}\bigr)_{l,l'<n}\cdot\Xi^{\T},
\end{align*}
and $\det\Xi=1$.
\end{proof}

\begin{corollary}\label{cor:cindep}
Let $A=\Z[x,c]$ and $s=(c+x,c,c,\dots)$, and write $J^{(c)},P^{(c)}_m,L^{(c)},a^{(c)}_{n,j}$ for the objects $J^{(s)},P^{(s)}_m,L^{(s)},a^{(s)}_{n,j}$ of Lemma~\ref{lem:jacobi}; thus $J^{(c)}$ is the matrix $J$ of Section~\ref{sec:prelim}, and $a^{(c)}_{n,j}=a_{n,j}(x,c)$ by Lemma~\ref{lem:transfer}. Then the matrix $P^{(c)}_k(J^{(c)})$ does not depend on $c$; hence $D_{k,n}(x,c)=D_{k,n}(x,0)$. $($This reproves the $c$-independence in Theorem~\ref{thm:A}.$)$
\end{corollary}

\begin{proof}
Let $J^{(0)},P^{(0)}_m$ be the objects for $c=0$, that is, for $s=(x,0,0,\dots)$. Then $J^{(c)}=J^{(0)}+cI$, and $P^{(c)}_m(X)=P^{(0)}_m(X-c)$ for all $m\ge0$ (induction on $m$ from the defining recurrences), so
$P^{(c)}_k(J^{(c)})=P^{(0)}_k(J^{(c)}-cI)=P^{(0)}_k(J^{(0)})$.
Now apply Lemma~\ref{lem:jacobi}(iv) to $s=(c+x,c,c,\dots)$ and to $s=(x,0,0,\dots)$.
\end{proof}

\emph{Conventions for the rest of this section.} By Corollary~\ref{cor:cindep} it suffices to prove Theorem~\ref{thm:B} for $c=0$. So from now on $c=0$, $s=(x,0,0,\dots)$, $K:=\Q(x)$ (the field of rational functions in the indeterminate $x$), $D_{k,n}=D_{k,n}(x,0)\in\Z[x]\subset K$, and $J,P_m,L$ (without superscript) denote the objects $J^{(s)},P^{(s)}_m,L^{(s)}$ of Lemma~\ref{lem:jacobi} for $A=K$ and $s=(x,0,0,\dots)$, that is, $J^{(0)},P^{(0)}_m,L^{(0)}$ regarded over $K$. The polynomials $P_m$ are
\begin{equation}\label{eq:Pm}
P_m=U_m-xU_{m-1}\qquad(m\ge0).
\end{equation}
Indeed, the right-hand side equals $1$ for $m=0$ and $X-x$ for $m=1$, and satisfies $P_{m+1}=XP_m-P_{m-1}$ for $m\ge1$ (because both $(U_m)$ and $(U_{m-1})$ satisfy this recurrence).

\subsection{A lemma on complementary minors}

\begin{lemma}[complementary minors]\label{lem:compl}
Let $r,r'\ge1$ be integers, $M\in K^{r\times(r+r')}$ and $Z\in K^{(r+r')\times r'}$ with $MZ=0$. Index rows and columns from $0$. For $S\subseteq\{0,\dots,r+r'-1\}$ with $|S|=r$ let $M_S$ be the $r\times r$ submatrix of $M$ formed by the columns with index in $S$, let $Z_S$ be the $r'\times r'$ submatrix of $Z$ formed by the rows with index not in $S$ $($rows and columns kept in increasing order$)$, and let $\|S\|$ denote the sum of the elements of $S$. Then for all such $S,S'$:
\[
(-1)^{\|S\|}\det(M_S)\det(Z_{S'})=(-1)^{\|S'\|}\det(M_{S'})\det(Z_S).
\]
\end{lemma}

This is the classical duality between the Plücker coordinates of the row space of $M$ and those of its annihilator; we include the short proof to fix the signs.

\begin{proof}
If $\rk M<r$ or $\rk Z<r'$, then all $\det M_S$ or all $\det Z_S$ vanish, and the identity reads $0=0$. Otherwise the transposes of the rows of $M$ form a basis of the kernel of the surjective map $\pi\colon K^{r+r'}\to K^{r'}$, $\pi(u)=Z^{\T}u$ (they lie in $\ker\pi$ because $Z^{\T}M^{\T}=(MZ)^{\T}=0$, they are linearly independent because $\rk M=r$, and $\dim\ker\pi=r$). Let $\psi(u_1,\dots,u_{r'})$ be the determinant of the square matrix whose first $r$ rows are the rows of $M$ and whose last $r'$ rows are $u_1^{\T},\dots,u_{r'}^{\T}$. The map $\psi$ is multilinear and alternating, and it vanishes whenever some $u_i$ lies in $\ker\pi$; hence it induces an alternating $r'$-linear form on $K^{r+r'}/\ker\pi\cong K^{r'}$, and there is $\psi_0\in K$ with $\psi(u_1,\dots,u_{r'})=\psi_0\det\bigl(\pi(u_1)|\cdots|\pi(u_{r'})\bigr)$. Let $S=\{i_0<\dots<i_{r-1}\}$, let $j_1<\dots<j_{r'}$ be the elements of its complement, and take $u_m=e_{j_m}$ $(1\le m\le r')$, where $e_0,\dots,e_{r+r'-1}$ is the standard basis of $K^{r+r'}$. Then $\det\bigl(\pi(u_1)|\cdots|\pi(u_{r'})\bigr)=\det Z_S$. On the other hand, moving the columns with index in $S$ to the front (keeping the order inside $S$ and inside its complement) turns the matrix defining $\psi(e_{j_1},\dots,e_{j_{r'}})$ into a block upper triangular matrix with diagonal blocks $M_S$ and $I_{r'}$, and this column permutation has sign $(-1)^{\sum_{m=0}^{r-1}(i_m-m)}=(-1)^{\|S\|-r(r-1)/2}$. Hence $(-1)^{\|S\|-r(r-1)/2}\det M_S=\psi_0\det Z_S$ for every $S$; multiplying this identity for $S$ by $\det Z_{S'}$ and comparing with the identity for $S'$ multiplied by $\det Z_S$ gives the claim.
\end{proof}

\subsection{Reduction of $D_{k,n}$ to a $k\times k$ determinant}
Let $R:=K[X]/(P_k)$, a $K$-vector space of dimension $k$. Since $P_0,\dots,P_{k-1}$ are monic of degrees $0,\dots,k-1$, their residue classes form a basis of $R$. For $f\in K[X]$ let $\crd(f)\in K^k$ (a row vector) be the coordinate vector of the residue class of $f$ with respect to this basis; $\crd$ is $K$-linear, $\crd(f)=0$ if $P_k$ divides $f$, and $\crd(P_0),\dots,\crd(P_{k-1})$ are the rows of the identity matrix $I_k$. For $f_1,\dots,f_k\in K[X]$ let $\detR(f_1,\dots,f_k)$ be the determinant of the $k\times k$ matrix with rows $\crd(f_1),\dots,\crd(f_k)$.

\begin{lemma}\label{lem:reduction}
For all $n\ge0$,
\[
D_{k,n}=(-1)^{nk}\detR(P_n,P_{n+1},\dots,P_{n+k-1}).
\]
\end{lemma}

Lemma~\ref{lem:reduction} is an instance of the Christoffel-type formulas for Hankel determinants of linear combinations of moments of orthogonal polynomials \cite{CK21,K23}; we give a direct proof.

\begin{proof}
For $n=0$ both sides are $1$. Let $n\ge1$ and $\Pi:=P_k(J)$, so that $\Pi_{j,m}=L(P_jP_kP_m)$ by Lemma~\ref{lem:jacobi}(ii). Fix $j\ge0$. The polynomial $P_kP_j$ is monic of degree $j+k$, and $(P_m)_{m\ge0}$ is a basis of $K[X]$ consisting of monic polynomials with $\deg P_m=m$, so $P_kP_j=\sum_{m=0}^{j+k}\nu_mP_m$ with $\nu_m\in K$ and $\nu_{j+k}=1$. Multiplying by $P_i$ and applying $L$, the orthogonality $L(P_iP_m)=\delta_{i,m}$ of Lemma~\ref{lem:jacobi}(ii) gives $\Pi_{j,i}=\nu_i$ for $i\le j+k$ and $\Pi_{j,i}=0$ for $i>j+k$. Thus
\[
\text{(a) } \Pi_{j,m}=0 \text{ for } m>j+k;\qquad \text{(b) } \Pi_{j,j+k}=1;\qquad \text{(c) } P_kP_j=\sum_{m=0}^{j+k}\Pi_{j,m}P_m.
\]
Applying $\crd$ to (c) gives $\sum_{m=0}^{j+k}\Pi_{j,m}\crd(P_m)=0$ for all $j\ge0$. Put $r=n$, $r'=k$,
\[
M=(\Pi_{j,m})_{0\le j\le n-1,\ 0\le m\le n+k-1}\in K^{n\times(n+k)},
\]
and let $Z$ be the $(n+k)\times k$ matrix with rows $\crd(P_0),\dots,\crd(P_{n+k-1})$. By (a) and the last identity, $MZ=0$ (for $j\le n-1$ all $m\le j+k$ are $\le n+k-1$). The hypothesis $r,r'\ge1$ of Lemma~\ref{lem:compl} holds because $n\ge1$ and $k\ge1$. Apply Lemma~\ref{lem:compl} with $S=\{0,\dots,n-1\}$ and $S'=\{k,\dots,n+k-1\}$:
\begin{itemize}
\item $\det M_S=\det(\Pi_{j,m})_{0\le j,m\le n-1}=D_{k,n}$ by Lemma~\ref{lem:jacobi}(iv);
\item $Z_S$ consists of the rows $\crd(P_n),\dots,\crd(P_{n+k-1})$, so $\det Z_S=\detR(P_n,\dots,P_{n+k-1})$;
\item $M_{S'}=(\Pi_{j,k+i})_{0\le j,i\le n-1}$ is lower unitriangular by (a) and (b), so $\det M_{S'}=1$;
\item $Z_{S'}$ consists of the rows $\crd(P_0),\dots,\crd(P_{k-1})$, that is, $Z_{S'}=I_k$ and $\det Z_{S'}=1$;
\item $\|S\|=n(n-1)/2$ and $\|S'\|=nk+n(n-1)/2$.
\end{itemize}
Lemma~\ref{lem:compl} gives $(-1)^{n(n-1)/2}D_{k,n}=(-1)^{nk+n(n-1)/2}\detR(P_n,\dots,P_{n+k-1})$.
\end{proof}

\subsection{The Laurent polynomial model}\label{ss:laurent}
Let $B:=K[t,t^{-1}]$ (the ring of Laurent polynomials), let $\sigma\colon B\to B$ be the $K$-algebra involution $t\mapsto t^{-1}$, and put
\begin{gather*}
B^+=\{f\in B:\sigma(f)=f\},\qquad B^-=\{f\in B:\sigma(f)=-f\},\\
y:=t+t^{-1}\in B^+,\qquad Q(t):=t^{2k+2}-x\,t^{2k+1}+x\,t-1\in B.
\end{gather*}
Recall \eqref{eq:cheb}: $(t-t^{-1})\,U_m(y)=t^{m+1}-t^{-m-1}$ for all $m\in\Z$.

\begin{lemma}\label{lem:laurent}
\begin{enumerate}
\item[(i)] $t^{-k-1}Q(t)=(t-t^{-1})\,P_k(y)$; in particular $\sigma(Q)=-t^{-2k-2}Q$.
\item[(ii)] The $K$-linear map $\mu\colon K[X]\to B^-$, $\mu(f)=(t-t^{-1})f(y)$, is bijective, and $\mu(P_kK[X])=QB\cap B^-$.
\item[(iii)] Let $V:=B^-/(QB\cap B^-)$, and for $j\in\Z$ let $F_j\in V$ be the class of $t^j-t^{-j}$. Then $\mu$ induces a $K$-linear isomorphism $\iota\colon R\to V$. Writing $\iota(f)$ for the image under $\iota$ of the residue class of $f\in K[X]$, we have $\iota(U_{j-1})=F_j$ for all $j\in\Z$ and $\iota(P_m)=F_{m+1}-xF_m$ for all $m\ge0$. Moreover $F_0=0$, $F_{-j}=-F_j$, and $(F_1,\dots,F_k)$ is a basis of $V$.
\item[(iv)] $F_{j+2k+2}-xF_{j+2k+1}+xF_{j+1}-F_j=0$ in $V$ for every $j\in\Z$.
\item[(v)] Let $\ell\in\Z$, $\gamma_0,\dots,\gamma_{k+1}\in K$ and $\gamma(t)=\sum_{i=0}^{k+1}\gamma_it^i$. Then $\sum_{i=0}^{k+1}\gamma_iF_{\ell-1+i}=0$ in $V$ if and only if $t^{\ell-1}\gamma(t)-t^{1-\ell}\gamma(t^{-1})\in QB$.
\end{enumerate}
\end{lemma}

\begin{proof}
(i) By \eqref{eq:cheb} and \eqref{eq:Pm}, $(t-t^{-1})P_k(y)=(t^{k+1}-t^{-k-1})-x(t^k-t^{-k})=t^{-k-1}Q(t)$. The left side is antisymmetric, so $\sigma(t^{-k-1}Q)=-t^{-k-1}Q$, that is, $\sigma(Q)=-t^{-2k-2}Q$.

(ii) The homomorphism $K[X]\to B$, $f\mapsto f(y)$, is injective (compare highest exponents) and $B$ is a domain, so $\mu$ is injective; its image lies in $B^-$ since $t-t^{-1}\in B^-$ and $f(y)\in B^+$. Surjectivity: let $f=\sum_jf_jt^j\in B^-$. Then $f_{-j}=-f_j$ for all $j$, so $f_0=0$ and $f=\sum_{j\ge1}f_j(t^j-t^{-j})=(t-t^{-1})\sum_{j\ge1}f_jU_{j-1}(y)$ by \eqref{eq:cheb}. Next, $B^+=\{f(y):f\in K[X]\}$: if $f=\sum_jf_jt^j\in B^+$ then $f=f_0+\sum_{j\ge1}f_j(t^j+t^{-j})$, and $t^j+t^{-j}=U_j(y)-U_{j-2}(y)$ for $j\ge1$ (multiply by $t-t^{-1}$ and use \eqref{eq:cheb}; recall $U_{-1}=0$). Finally, since $t^{k+1}$ is a unit of $B$, (i) gives $QB=(t-t^{-1})P_k(y)B$. Let $f\in B$. Since $\sigma$ is a ring homomorphism with $\sigma(t-t^{-1})=-(t-t^{-1})$ and $\sigma(y)=y$, we have
\[
\sigma\bigl((t-t^{-1})P_k(y)f\bigr)=-(t-t^{-1})P_k(y)\,\sigma(f).
\]
Hence $(t-t^{-1})P_k(y)f$ lies in $B^-$ if and only if $(t-t^{-1})P_k(y)\bigl(\sigma(f)-f\bigr)=0$, that is, if and only if $\sigma(f)=f$ ($B$ is a domain and $(t-t^{-1})P_k(y)=\mu(P_k)\ne0$). Hence $QB\cap B^-=(t-t^{-1})P_k(y)B^+=\mu(P_kK[X])$.

(iii) By (ii), $\mu$ induces an isomorphism $\iota\colon K[X]/(P_k)\to B^-/(QB\cap B^-)$. By \eqref{eq:cheb}, $\mu(U_{j-1})=t^j-t^{-j}$ for all $j\in\Z$, so $\iota(U_{j-1})=F_j$, and by \eqref{eq:Pm} $\iota(P_m)=F_{m+1}-xF_m$ for $m\ge0$. The relations $F_0=0$ and $F_{-j}=-F_j$ are clear from the definition. The residue classes of $U_0,\dots,U_{k-1}$ (monic of degrees $0,\dots,k-1$) form a basis of $R$, and $\iota$ maps them to $F_1,\dots,F_k$.

(iv) We have $(t^j+t^{-j-2k-2})Q(t)\in QB$, and expanding,
\begin{align*}
(t^j+t^{-j-2k-2})Q(t)={}&\bigl(t^{j+2k+2}-t^{-(j+2k+2)}\bigr)-x\bigl(t^{j+2k+1}-t^{-(j+2k+1)}\bigr)\\
&+x\bigl(t^{j+1}-t^{-(j+1)}\bigr)-\bigl(t^j-t^{-j}\bigr),
\end{align*}
which lies in $B^-$. So it lies in $QB\cap B^-$, and its class in $V$, which is the left-hand side of (iv), is $0$.

(v) The element $\sum_i\gamma_i\bigl(t^{\ell-1+i}-t^{-(\ell-1+i)}\bigr)=t^{\ell-1}\gamma(t)-t^{1-\ell}\gamma(t^{-1})$ lies in $B^-$; its class in $V$ is $\sum_i\gamma_iF_{\ell-1+i}$, and by the definition of $V$ this class is $0$ if and only if the element lies in $QB$.
\end{proof}

For $u_1,\dots,u_k\in V$ let $\detV(u_1,\dots,u_k)$ be the determinant of the $k\times k$ matrix whose $i$-th row is the coordinate vector of $u_i$ with respect to the basis $(F_1,\dots,F_k)$ of $V$. For $m\in\Z$ put
\begin{equation}\label{eq:wG}
w(m):=\detV(F_m,F_{m+1},\dots,F_{m+k-1}),\qquad G_m:=F_{m+1}-xF_m\in V.
\end{equation}
By Lemma~\ref{lem:laurent}(iii), $w(m)$ is exactly the quantity $w(m)$ of Theorem~\ref{thm:C}: $\iota$ maps the residue class of $U_{j-1}$ to $F_j$ for every $j\in\Z$, and it maps the basis of $R$ formed by the residue classes of $U_0,\dots,U_{k-1}$ to the basis $(F_1,\dots,F_k)$ of $V$; so the coordinate vector of the residue class of $U_{j-1}$ modulo $P_k$ with respect to the residue classes of $U_0,\dots,U_{k-1}$ equals the coordinate vector of $F_j$ with respect to $(F_1,\dots,F_k)$.

\begin{lemma}\label{lem:Lambda}
\begin{enumerate}
\item[(i)] $\detR(f_1,\dots,f_k)=\detV\bigl(\iota(f_1),\dots,\iota(f_k)\bigr)$ for all $f_1,\dots,f_k\in K[X]$. Consequently
\[
D_{k,n}=(-1)^{nk}\detV(G_n,G_{n+1},\dots,G_{n+k-1})\qquad(n\ge0).
\]
\item[(ii)] There is a $K$-linear map $\Lambda\colon V\to V$ such that for all $j\in\Z$
\[
(\Lambda1)\quad \Lambda F_j=F_{j+k+1}-xF_{j+k},\qquad\qquad (\Lambda2)\quad \Lambda F_j=F_{j-k-1}-xF_{j-k}.
\]
\item[(iii)] $\det\Lambda=(-1)^{k(k+1)/2}$; in particular $(\det\Lambda)^2=1$.
\item[(iv)] $D_{k,n}=(-1)^{nk+k(k+1)/2}\,w(n-k)$ for all $n\ge0$.
\end{enumerate}
\end{lemma}

\begin{proof}
(i) The coordinates of $\iota(f)$ with respect to the basis $(\iota(P_0),\dots,\iota(P_{k-1}))$ of $V$ are $\crd(f)$. By Lemma~\ref{lem:laurent}(iii), $\iota(P_0)=F_1$ (as $F_0=0$) and $\iota(P_m)=F_{m+1}-xF_m$ for $1\le m\le k-1$, so the basis $(\iota(P_m))_{m<k}$ is obtained from $(F_1,\dots,F_k)$ by a unitriangular matrix, whose determinant is $1$. Hence the two determinants agree. The formula for $D_{k,n}$ follows from Lemma~\ref{lem:reduction} and $\iota(P_m)=G_m$ $(m\ge0)$.

(ii) Let $\vartheta:=t^{k+1}-xt^k+t^{-k-1}-xt^{-k}\in B^+$. Multiplication by $\vartheta$ maps $B^-$ into $B^-$ and $QB\cap B^-$ into itself, so it induces a $K$-linear map $\vartheta_V\colon V\to V$. Expanding the product,
\[
\vartheta\cdot(t^j-t^{-j})=\Bigl[\bigl(t^{j+k+1}-t^{-(j+k+1)}\bigr)-x\bigl(t^{j+k}-t^{-(j+k)}\bigr)\Bigr]+\Bigl[\bigl(t^{j-k-1}-t^{-(j-k-1)}\bigr)-x\bigl(t^{j-k}-t^{-(j-k)}\bigr)\Bigr],
\]
so $\vartheta_VF_j=(F_{j+k+1}-xF_{j+k})+(F_{j-k-1}-xF_{j-k})$. By Lemma~\ref{lem:laurent}(iv), applied with $j-k-1$ in place of $j$, the two brackets are equal in $V$: $F_{j+k+1}-xF_{j+k}=F_{j-k-1}-xF_{j-k}$. Hence $\Lambda:=\frac12\vartheta_V$ satisfies ($\Lambda1$) and ($\Lambda2$).

(iii) By ($\Lambda1$) with $j=m-k$, $\Lambda F_{m-k}=F_{m+1}-xF_m=G_m$ for all $m\in\Z$. The vectors $G_0=F_1$, $G_1=F_2-xF_1$, \dots, $G_{k-1}=F_k-xF_{k-1}$ have a unitriangular coordinate matrix, so
\begin{align*}
1&=\detV(G_0,\dots,G_{k-1})=\detV(\Lambda F_{-k},\dots,\Lambda F_{-1})=\det(\Lambda)\cdot\detV(F_{-k},\dots,F_{-1})\\
&=\det(\Lambda)\cdot\detV(-F_k,-F_{k-1},\dots,-F_1)=\det(\Lambda)\cdot(-1)^k(-1)^{k(k-1)/2}\detV(F_1,\dots,F_k)\\
&=\det(\Lambda)\cdot(-1)^{k(k+1)/2}.
\end{align*}
(Here $\detV(\Lambda u_1,\dots,\Lambda u_k)=\det(\Lambda)\detV(u_1,\dots,u_k)$ for all $u_1,\dots,u_k\in V$, and reversing the order of $k$ rows has sign $(-1)^{k(k-1)/2}$.)

(iv) By (i) and $G_m=\Lambda F_{m-k}$,
\[
D_{k,n}=(-1)^{nk}\detV(\Lambda F_{n-k},\dots,\Lambda F_{n-1})=(-1)^{nk}\det(\Lambda)\,w(n-k);
\]
now use (iii).
\end{proof}

\subsection{The bilinear identity}

\begin{proposition}\label{prop:bilinear}
For every $\ell\in\Z$,
\[
w(\ell+k+1)\,w(\ell-k-1)-x^2\,w(\ell+k)\,w(\ell-k)=(1-x^2)\,w(\ell)^2.
\]
\end{proposition}

\begin{proof}
Fix $\ell\in\Z$ and put $v_i:=F_{\ell-1+i}\in V$ for $i=0,1,\dots,k+1$ (a ``window'' of $k+2$ consecutive vectors). Let $C$ be the $k\times(k+2)$ matrix whose $i$-th column is the coordinate vector of $v_i$ with respect to $(F_1,\dots,F_k)$, and for $0\le i<j\le k+1$ let $[i,j]$ be the determinant of the $k\times k$ matrix obtained from $C$ by deleting columns $i$ and $j$. (Since $\detV$ is defined via rows and $C$ via columns, $\detV$ of a subfamily of the $v_i$ in increasing order is the corresponding minor of $C$.)

\emph{Step 1 (the five determinants as combinations of the minors $[i,j]$).} Define
\begin{equation}\label{eq:Deltas}
\begin{aligned}
\Delta_1&:=\detV(v_1-xv_2,\ v_2-xv_3,\ \dots,\ v_k-xv_{k+1}),\\
\Delta_2&:=\detV(v_1-xv_0,\ v_2-xv_1,\ \dots,\ v_k-xv_{k-1}),\\
\Delta_3&:=\detV(v_0-xv_1,\ v_1-xv_2,\ \dots,\ v_{k-1}-xv_k),\\
\Delta_4&:=\detV(v_2-xv_1,\ v_3-xv_2,\ \dots,\ v_{k+1}-xv_k),\\
\Delta_0&:=\detV(v_1,v_2,\dots,v_k)=w(\ell)=[0,k+1].
\end{aligned}
\end{equation}

\emph{Claim 1a.} We have
\begin{equation}\label{eq:claim1a}
\begin{aligned}
\Delta_1&=\det(\Lambda)\,w(\ell+k+1), &\qquad \Delta_2&=\det(\Lambda)\,w(\ell-k-1),\\
\Delta_3&=\det(\Lambda)\,w(\ell+k), &\qquad \Delta_4&=\det(\Lambda)\,w(\ell-k).
\end{aligned}
\end{equation}
Indeed, for $i=0,\dots,k-1$:
\begin{itemize}
\item by ($\Lambda2$), $\Lambda F_{\ell+k+1+i}=F_{\ell+i}-xF_{\ell+i+1}=v_{i+1}-xv_{i+2}$;
\item by ($\Lambda1$), $\Lambda F_{\ell-k-1+i}=F_{\ell+i}-xF_{\ell+i-1}=v_{i+1}-xv_i$;
\item by ($\Lambda2$), $\Lambda F_{\ell+k+i}=F_{\ell+i-1}-xF_{\ell+i}=v_i-xv_{i+1}$;
\item by ($\Lambda1$), $\Lambda F_{\ell-k+i}=F_{\ell+i+1}-xF_{\ell+i}=v_{i+2}-xv_{i+1}$;
\end{itemize}
and $\detV(\Lambda u_1,\dots,\Lambda u_k)=\det(\Lambda)\detV(u_1,\dots,u_k)$.

\emph{Claim 1b.} We have
\begin{equation}\label{eq:claim1b}
\begin{aligned}
\Delta_1&=\sum_{m=1}^{k+1}(-x)^{k+1-m}\,[0,m], &\qquad \Delta_2&=\sum_{m=0}^{k}(-x)^m\,[m,k+1],\\
\Delta_3&=\sum_{m=0}^{k}(-x)^{k-m}\,[m,k+1], &\qquad \Delta_4&=\sum_{m=1}^{k+1}(-x)^{m-1}\,[0,m].
\end{aligned}
\end{equation}
To prove Claim 1b, expand each determinant by multilinearity into $2^k$ terms, one for each choice, in every argument, of its first or its second summand; a term vanishes if two of the chosen vectors coincide.
\begin{itemize}
\item $\Delta_1$: argument $j$ $(j=1,\dots,k)$ is $v_j-xv_{j+1}$. The chosen indices are pairwise distinct if and only if there is $m\in\{1,\dots,k+1\}$ such that the arguments $j<m$ choose $v_j$ and the arguments $j\ge m$ choose $-xv_{j+1}$. The chosen vectors are then $v_1,\dots,v_{m-1},v_{m+1},\dots,v_{k+1}$ in increasing order, giving $(-x)^{k+1-m}[0,m]$.
\item $\Delta_2$: argument $j$ $(j=1,\dots,k)$ is $v_j-xv_{j-1}$. The chosen indices are pairwise distinct if and only if there is $m\in\{0,\dots,k\}$ such that the arguments $j\le m$ choose $-xv_{j-1}$ and the arguments $j>m$ choose $v_j$: the vectors are $v_0,\dots,v_{m-1},v_{m+1},\dots,v_k$, giving $(-x)^m[m,k+1]$.
\item $\Delta_3$: argument $j$ $(j=0,\dots,k-1)$ is $v_j-xv_{j+1}$. The chosen indices are pairwise distinct if and only if there is $m\in\{0,\dots,k\}$ such that the arguments $j<m$ choose $v_j$ and the arguments $j\ge m$ choose $-xv_{j+1}$: the vectors are $v_0,\dots,v_{m-1},v_{m+1},\dots,v_k$, giving $(-x)^{k-m}[m,k+1]$.
\item $\Delta_4$: argument $j$ $(j=1,\dots,k)$ is $v_{j+1}-xv_j$. The chosen indices are pairwise distinct if and only if there is $m\in\{1,\dots,k+1\}$ such that the arguments $j<m$ choose $-xv_j$ and the arguments $j\ge m$ choose $v_{j+1}$: the vectors are $v_1,\dots,v_{m-1},v_{m+1},\dots,v_{k+1}$, giving $(-x)^{m-1}[0,m]$.
\end{itemize}
(In each case, if the choices are not of the stated ``monotone'' form, two consecutive arguments choose the same vector.)

By Claim 1a and $(\det\Lambda)^2=1$ (Lemma~\ref{lem:Lambda}(iii)),
\begin{equation}\label{eq:reduce}
w(\ell+k+1)\,w(\ell-k-1)-x^2\,w(\ell+k)\,w(\ell-k)=\Delta_1\Delta_2-x^2\Delta_3\Delta_4,
\end{equation}
so it remains to prove
\begin{equation}\label{eq:target}
\Delta_1\Delta_2-x^2\Delta_3\Delta_4=(1-x^2)\,\Delta_0^2.
\end{equation}

\emph{Step 2 (the kernel of the window).} If $\rk C<k$, then all minors $[i,j]$ vanish, so $\Delta_1=\Delta_2=\Delta_3=\Delta_4=\Delta_0=0$ by Claim 1b and \eqref{eq:Deltas}, and \eqref{eq:target} holds. (Such degenerate windows do occur; see Remark~\ref{rem:degenerate}.) Assume from now on that $\rk C=k$. Then the kernel of $C$ in $K^{k+2}$ has dimension $2$; choose a basis $\gamma=(\gamma_0,\dots,\gamma_{k+1})^{\T}$, $\eta=(\eta_0,\dots,\eta_{k+1})^{\T}$ of it, and put $g_i:=(\gamma_i,\eta_i)^{\T}\in K^2$ $(0\le i\le k+1)$. For $u_1,u_2\in K^2$ write $\det(u_1,u_2)$ for the determinant of the $2\times2$ matrix with columns $u_1,u_2$; thus $\det(g_i,g_j)=\gamma_i\eta_j-\gamma_j\eta_i$. Apply Lemma~\ref{lem:compl} with $r=k$, $r'=2$, $M=C$ and $Z$ the $(k+2)\times2$ matrix with columns $\gamma,\eta$ (so $CZ=0$). For $S=\{0,\dots,k+1\}\setminus\{i,j\}$ $(i<j)$ we get $\det M_S=[i,j]$, $\det Z_S=\det(g_i,g_j)$ and $(-1)^{\|S\|}=(-1)^{(k+1)(k+2)/2}(-1)^{i+j}$. Hence for all pairs $i<j$ and $i'<j'$,
\begin{equation}\label{eq:pluck-ratio}
(-1)^{i+j}\,[i,j]\,\det(g_{i'},g_{j'})=(-1)^{i'+j'}\,[i',j']\,\det(g_i,g_j).
\end{equation}
Since $\gamma,\eta$ are linearly independent, $\det(g_{i'},g_{j'})\ne0$ for some pair $i'<j'$; with
\[
\kappa:=\frac{(-1)^{i'+j'}[i',j']}{\det(g_{i'},g_{j'})}
\]
we obtain
\begin{equation}\label{eq:kappa}
[i,j]=\kappa\,(-1)^{i+j}\det(g_i,g_j)\qquad\text{for all }0\le i<j\le k+1.
\end{equation}

\emph{Step 3 (rewriting \eqref{eq:target}).} Let $\gamma(t)=\sum_{i=0}^{k+1}\gamma_it^i$ and $\eta(t)=\sum_{i=0}^{k+1}\eta_it^i$ (polynomials in $t$ of degree at most $k+1$, elements of $B$) and, for $\tau\in K$, $g(\tau):=(\gamma(\tau),\eta(\tau))^{\T}=\sum_{i=0}^{k+1}\tau^ig_i\in K^2$. Recall that $x\ne0$ in $K$. Inserting \eqref{eq:kappa} into Claim 1b, and using
\begin{gather*}
(-x)^{k+1-m}(-1)^m=(-1)^{k+1}x^{k+1-m},\qquad (-x)^m(-1)^{m+k+1}=(-1)^{k+1}x^m,\\
(-x)^{k-m}(-1)^{m+k+1}=-x^{k-m},\qquad (-x)^{m-1}(-1)^m=-x^{m-1},
\end{gather*}
the bilinearity of $\det(\cdot,\cdot)$ and $\det(u_1,u_1)=0$, we obtain
\begin{equation}\label{eq:Deltas-g}
\begin{aligned}
\Delta_1&=\kappa(-1)^{k+1}\det\Bigl(g_0,\sum_{m=1}^{k+1}x^{k+1-m}g_m\Bigr)=\kappa(-1)^{k+1}x^{k+1}\det\bigl(g_0,g(1/x)\bigr),\\
\Delta_2&=\kappa(-1)^{k+1}\det\Bigl(\sum_{m=0}^{k}x^mg_m,\,g_{k+1}\Bigr)=\kappa(-1)^{k+1}\det\bigl(g(x),g_{k+1}\bigr),\\
\Delta_3&=-\kappa\det\Bigl(\sum_{m=0}^{k}x^{k-m}g_m,\,g_{k+1}\Bigr)=-\kappa\,x^k\det\bigl(g(1/x),g_{k+1}\bigr),\\
\Delta_4&=-\kappa\det\Bigl(g_0,\sum_{m=1}^{k+1}x^{m-1}g_m\Bigr)=-\kappa\,x^{-1}\det\bigl(g_0,g(x)\bigr),\\
\Delta_0&=\kappa(-1)^{k+1}\det(g_0,g_{k+1}).
\end{aligned}
\end{equation}
(Here $\sum_{m=1}^{k+1}x^{k+1-m}g_m=x^{k+1}\bigl(g(1/x)-g_0\bigr)$, $\sum_{m=0}^kx^mg_m=g(x)-x^{k+1}g_{k+1}$, $\sum_{m=0}^kx^{k-m}g_m=x^k\bigl(g(1/x)-x^{-k-1}g_{k+1}\bigr)$ and $\sum_{m=1}^{k+1}x^{m-1}g_m=x^{-1}\bigl(g(x)-g_0\bigr)$.) Therefore
\begin{equation}\label{eq:bracket}
\Delta_1\Delta_2-x^2\Delta_3\Delta_4=\kappa^2x^{k+1}\Bigl[\det\bigl(g_0,g(1/x)\bigr)\det\bigl(g(x),g_{k+1}\bigr)-\det\bigl(g_0,g(x)\bigr)\det\bigl(g(1/x),g_{k+1}\bigr)\Bigr].
\end{equation}
The three-term Plücker relation for four vectors $u_1,u_2,u_3,u_4\in K^2$,
\[
\det(u_1,u_2)\det(u_3,u_4)-\det(u_1,u_3)\det(u_2,u_4)+\det(u_1,u_4)\det(u_2,u_3)=0
\]
(direct verification), applied to $u_1=g_0$, $u_2=g(1/x)$, $u_3=g(x)$, $u_4=g_{k+1}$, turns the bracket into $-\det(g_0,g_{k+1})\det\bigl(g(1/x),g(x)\bigr)=\det(g_0,g_{k+1})\det\bigl(g(x),g(1/x)\bigr)$. Hence
\begin{equation}\label{eq:after-pluecker}
\Delta_1\Delta_2-x^2\Delta_3\Delta_4=\kappa^2x^{k+1}\det(g_0,g_{k+1})\cdot\det\bigl(g(x),g(1/x)\bigr).
\end{equation}

\emph{Step 4 (the key divisibility).} Consider the Laurent polynomial
\[
\Phi(t):=\gamma(t)\,\eta(t^{-1})-\eta(t)\,\gamma(t^{-1})\in B,
\]
all of whose exponents lie in $[-(k+1),k+1]$. Since $C\gamma=0$ and $C\eta=0$, that is, $\sum_i\gamma_iv_i=\sum_i\eta_iv_i=0$ in $V$, Lemma~\ref{lem:laurent}(v) gives $\alpha,\beta\in B$ with
\begin{equation}\label{eq:alphabeta}
t^{\ell-1}\gamma(t)-t^{1-\ell}\gamma(t^{-1})=Q\alpha,\qquad t^{\ell-1}\eta(t)-t^{1-\ell}\eta(t^{-1})=Q\beta.
\end{equation}
Thus $\gamma(t^{-1})=t^{2\ell-2}\gamma(t)-t^{\ell-1}Q\alpha$ and $\eta(t^{-1})=t^{2\ell-2}\eta(t)-t^{\ell-1}Q\beta$, and substituting,
\begin{equation}\label{eq:Phi}
\Phi=\gamma(t)\bigl[t^{2\ell-2}\eta(t)-t^{\ell-1}Q\beta\bigr]-\eta(t)\bigl[t^{2\ell-2}\gamma(t)-t^{\ell-1}Q\alpha\bigr]=t^{\ell-1}Q(t)\bigl(\eta(t)\alpha-\gamma(t)\beta\bigr)\in QB.
\end{equation}
Write $\Phi=Q\zeta$ with $\zeta\in B$. If $\zeta\ne0$, let $m_1\le m_2$ be the lowest and the highest exponent occurring in $\zeta$. Since $Q$ has constant term $-1$ and leading term $t^{2k+2}$, the lowest and highest exponents of $Q\zeta$ are $m_1$ and $m_2+2k+2$. As the exponents of $\Phi$ lie in $[-(k+1),k+1]$, we get $m_1\ge-(k+1)$ and $m_2\le-(k+1)$, so $m_1=m_2=-(k+1)$. In all cases $\Phi=\lambda\,t^{-k-1}Q(t)$ for some $\lambda\in K$. Comparing the coefficients of $t^{k+1}$ (on the left: $\gamma_{k+1}\eta_0-\eta_{k+1}\gamma_0$; on the right: $\lambda$) gives $\lambda=-\det(g_0,g_{k+1})$, that is,
\begin{equation}\label{eq:Phi-final}
\gamma(t)\,\eta(t^{-1})-\eta(t)\,\gamma(t^{-1})=-\det(g_0,g_{k+1})\cdot t^{-k-1}Q(t).
\end{equation}
Now apply the $K$-algebra homomorphism $B\to K$, $t\mapsto x$ (well defined because $x\ne0$). The left side of \eqref{eq:Phi-final} becomes $\gamma(x)\eta(1/x)-\eta(x)\gamma(1/x)=\det\bigl(g(x),g(1/x)\bigr)$, and $Q(x)=x^{2k+2}-x^{2k+2}+x^2-1=x^2-1$. So
\begin{equation}\label{eq:at-x}
\det\bigl(g(x),g(1/x)\bigr)=(1-x^2)\,x^{-k-1}\det(g_0,g_{k+1}).
\end{equation}

\emph{Step 5.} Inserting \eqref{eq:at-x} into \eqref{eq:after-pluecker} gives $\Delta_1\Delta_2-x^2\Delta_3\Delta_4=\kappa^2(1-x^2)\det(g_0,g_{k+1})^2=(1-x^2)\Delta_0^2$ by the formula for $\Delta_0$ in \eqref{eq:Deltas-g}. This is \eqref{eq:target}.
\end{proof}

\subsection{Completion of the proofs}

\begin{proof}[Proof of Theorem~\ref{thm:C}]
The first formula is Lemma~\ref{lem:Lambda}(iv) (which is stated for $D_{k,n}(x,0)$; by Corollary~\ref{cor:cindep} this equals $D_{k,n}(x,c)$), and the second formula is Proposition~\ref{prop:bilinear}.
\end{proof}

Both identities of Theorem~\ref{thm:C} take place in $\Z[x]$, because $w(m)\in\Z[x]$ for every $m\in\Z$. Indeed, $P_k=U_k-xU_{k-1}$ is monic in $X$ with coefficients in $\Z[x]$, so division with remainder by $P_k$ can be carried out in $\Z[x][X]$, and the residue class of each $U_{j-1}\in\Z[X]$ $(j\in\Z)$ modulo $P_k$ is represented by a polynomial in $\Z[x][X]$ of degree less than $k$. Since $U_0,\dots,U_{k-1}$ are monic of degrees $0,\dots,k-1$ with coefficients in $\Z$, the coordinates of this residue class with respect to the residue classes of $U_0,\dots,U_{k-1}$ lie in $\Z[x]$. Thus $w(m)$ is the determinant of a matrix with entries in $\Z[x]$.

\begin{proof}[Proof of Theorem~\ref{thm:B}]
The case $k=0$ was treated at the beginning of this section; let $k\ge1$. By Corollary~\ref{cor:cindep} we may take $c=0$. Let $n\ge2k+2$ and $\ell:=n-2k-1\ge1$. All of $n$, $n-1$, $n-k-1$, $n-2k-1$, $n-2k-2$ are $\ge0$, so Lemma~\ref{lem:Lambda}(iv) applies to each of them; with $\varepsilon:=(-1)^{k(k+1)/2}$,
\begin{align*}
D_{k,n}&=(-1)^{nk}\varepsilon\,w(\ell+k+1), & D_{k,n-2k-2}&=(-1)^{(n-2k-2)k}\varepsilon\,w(\ell-k-1),\\
D_{k,n-1}&=(-1)^{(n-1)k}\varepsilon\,w(\ell+k), & D_{k,n-2k-1}&=(-1)^{(n-2k-1)k}\varepsilon\,w(\ell-k),\\
D_{k,n-k-1}&=(-1)^{(n-k-1)k}\varepsilon\,w(\ell).
\end{align*}
The sign of $D_{k,n}D_{k,n-2k-2}$ and the sign of $D_{k,n-1}D_{k,n-2k-1}$ are both $(-1)^{2nk-2k^2-2k}\varepsilon^2=1$, and $(D_{k,n-k-1})^2=w(\ell)^2$. Hence, by Proposition~\ref{prop:bilinear},
\begin{align*}
D_{k,n}D_{k,n-2k-2}-x^2D_{k,n-1}D_{k,n-2k-1}&=w(\ell+k+1)\,w(\ell-k-1)-x^2\,w(\ell+k)\,w(\ell-k)\\
&=(1-x^2)\,w(\ell)^2=(1-x^2)(D_{k,n-k-1})^2.
\end{align*}
This is an identity in $K=\Q(x)$ between elements of $\Z[x]$ (the $D_{k,n}(x,0)$ lie in $\Z[x]$ by definition, and so do the $w(m)$, as shown above); hence it holds in $\Z[x]$, and by Corollary~\ref{cor:cindep} in $\Z[x,c]$ for $D_{k,n}(x,c)$. By Theorem~\ref{thm:A} it also holds for $d_{k,n}(x)=\det H_{k,n}(x)$.
\end{proof}

\section{Remarks}\label{sec:remarks}

In this section $D_{k,n}=D_{k,n}(x,c)$, which does not depend on $c$ by Corollary~\ref{cor:cindep}, and $K$, $R$, $P_k$ and $w(m)$ are as in Section~\ref{sec:proofBC}; $U_m(x)\in\Z[x]$ denotes the value of $U_m$ at $X=x$.

\begin{remark}[The cases $k=1,2$; linear recurrences]\label{rem:smallk}
The question \cite{CMO} also displays closed formulas for $D_{1,n}$ and $D_{2,n}$; both follow from Theorem~\ref{thm:C}. For $k=1$ we have $P_1=X-x$, the residue class of $f\in K[X]$ modulo $P_1$ is $f(x)$ times the class of $U_0=1$, so $w(m)=U_{m-1}(x)$ and
\[
D_{1,n}=(-1)^{n+1}U_{n-2}(x)\qquad(n\ge0).
\]
Together with $U_n=\sum_{0\le2j\le n}(-1)^j\binom{n-j}{j}X^{n-2j}$ this gives Cigler's formula
\[
D_{1,n+2}=\sum_{0\le2j\le n}(-1)^{n-j-1}\binom{n-j}{j}x^{n-2j}\qquad(n\ge0),
\]
and $\sum_{n\ge0}D_{1,n}z^n=(1+xz)/(1+xz+z^2)$.

For $k=2$ we have $P_2=X^2-xX-1$ and $D_{2,n}=-w(n-2)$. Let $\rho_1,\rho_2$ be the roots of $P_2$ and choose $\omega_i$ with $\omega_i+\omega_i^{-1}=\rho_i$; here $\rho_1\ne\rho_2$ and $\omega_i\ne\pm1$ (see Remark~\ref{rem:symplectic}). Evaluating the residue classes at $\rho_1,\rho_2$ and using \eqref{eq:cheb} with $t\mapsto\omega_i$ shows that $w(m)$ is a linear combination, with coefficients independent of $m$, of the sequences $(\chi^m)_{m\in\Z}$ with $\chi\in\{\omega_1^{\pm1}\omega_2^{\pm1}\}$; since $\rho_1+\rho_2=x$ and $\rho_1\rho_2=-1$, each such $\chi$ satisfies $\chi^4+\chi^3+x^2\chi^2+\chi+1=0$. Hence
\begin{gather*}
D_{2,n+4}+D_{2,n+3}+x^2D_{2,n+2}+D_{2,n+1}+D_{2,n}=0\quad(n\ge0),\\
(D_{2,0},D_{2,1},D_{2,2},D_{2,3})=(1,0,0,-1),
\end{gather*}
so that $\sum_{n\ge0}D_{2,n}z^n=(1+z+x^2z^2)/(1+z+x^2z^2+z^3+z^4)$. Writing $1+z+x^2z^2+z^3+z^4=(1+z)(1+z^3)+x^2z^2$ and expanding $\sum_{n\ge0}D_{2,n+3}z^n=-(1+z)/\bigl((1+z)(1+z^3)+x^2z^2\bigr)$ as a geometric series in $x^2z^2(1+z)^{-1}(1+z^3)^{-1}$ gives Cigler's second formula
\[
D_{2,n+3}=\sum_{0\le2j\le n}(-1)^{n-j-1}\Bigl([z^{n-2j}](1-z)^{-j}(1-z^3)^{-j-1}\Bigr)x^{2j}\qquad(n\ge0).
\]
For $x=-t$ the two generating functions above are (2.5) and (2.6) of \cite{C22}, whose power series variable, called $x$ there, is our $z$; the quantity $d^{(k)}_0(n,t,s)$ of \cite{C22}, with level weight $s$ at height $0$ and $t$ at heights $>0$, is $D_{k,n}(s-t,t)$ in our notation.

For general $k$, Lemma~\ref{lem:jacobi}(iii) shows that $a_{n,k}=L^{(c)}(X^nP^{(c)}_k)$ is a linear combination, with leading coefficient $1$, of the $k+1$ consecutive moments $L^{(c)}(X^n),\dots,L^{(c)}(X^{n+k})$. This is the situation of \cite[Corollary~9]{K23} with $d=k$, $s=c$ and $t=1$; the scaling factor there is $\det\bigl(L^{(c)}(X^{i+j})\bigr)_{0\le i,j\le n-1}=1$ (Lemma~\ref{lem:jacobi}(iv) with $k=0$), and since all $t_i$ are equal the recurrence holds for all $n\ge0$ \cite[Remarks (1), (2) after Corollary~9]{K23}. Hence $(D_{k,n})_{n\ge0}$ satisfies a linear recurrence $\sum_{i=0}^{2^k}c_iD_{k,n+2^k-i}=0$ with constant coefficients, whose characteristic polynomial $\sum_{i=0}^{2^k}c_iY^{2^k-i}$ is that of the tensor product of the $k$ matrices $\bigl(\begin{smallmatrix}-\rho_j&1\\-1&0\end{smallmatrix}\bigr)$, $1\le j\le k$, where $\rho_1,\dots,\rho_k$ are the roots of $P_k$ (the roots of $P^{(c)}_k(X)=P_k(X-c)$ are $\rho_j+c$). For $k=1,2$ this gives the first two lines of (17) of \cite{C23}; for $k=3$ an exact computation gives the coefficients $1,-x,0,-(x^3-x),x^4+x^2-2,-(x^3-x),0,-x,1$, which are those of the third line.
\end{remark}

\begin{remark}[Degenerate windows]\label{rem:degenerate}
The identity of Proposition~\ref{prop:bilinear} holds for every integer $\ell$. Windows of rank less than $k$ (degenerate windows; this is the first case of Step~2 of the proof of Proposition~\ref{prop:bilinear}, in which all five determinants $\Delta_0,\dots,\Delta_4$ vanish) do occur, and not only near $\ell=0$:
\begin{itemize}
\item For $k\ge3$ and $-k+2\le\ell\le-1$ the window $F_{\ell-1},\dots,F_{\ell+k}$ is degenerate: all its indices lie in $[-(k-1),k-1]$, and $F_0=0$, $F_{-j}=-F_j$, so the window is contained in the span of $F_1,\dots,F_{k-1}$.
\item Exact computation with $x$ an indeterminate shows the following. For $k=1,2$ no window with $-4k-1\le\ell\le6k+1$ is degenerate; for $k=3$ the only degenerate window in this range is $\ell=-1$; for each $4\le k\le10$ there are also degenerate windows with $\ell\ge1$, that is, inside the range $\ell=n-2k-1\ge1$ which is used for Theorem~\ref{thm:B}, the smallest being $\ell=3$ (that is, $n=2k+4$). For example, for $k=4$ the degenerate windows with $-17\le\ell\le25$ are $\ell=-6,-2,-1,3$, and for $k=10$ exactly $53$ of the $103$ windows with $-41\le\ell\le61$ are degenerate.
\end{itemize}
Theorem~\ref{thm:B} needs $n\ge2k+2$ only so that all five indices $n$, $n-1$, $n-k-1$, $n-2k-1$, $n-2k-2$ are non-negative, which is needed to apply Lemma~\ref{lem:Lambda}(iv) to each of them.
\end{remark}

\begin{remark}[The operator $\Lambda$]\label{rem:Lambda}
The operator $\Lambda=\frac12\vartheta_V$ of Lemma~\ref{lem:Lambda}(ii) can be described in two further ways.
\begin{itemize}
\item Since $t^j+t^{-j}=U_j(y)-U_{j-2}(y)$ for $j\ge1$ (proof of Lemma~\ref{lem:laurent}(ii)), \eqref{eq:Pm} gives $\vartheta=(U_{k+1}-U_{k-1})(y)-x(U_k-U_{k-2})(y)=(P_{k+1}-P_{k-1})(y)$. As $\mu(Xf)=y\,\mu(f)$, the isomorphism $\iota$ satisfies $\iota(f_1f_2)=f_1(y)\iota(f_2)$ for $f_1,f_2\in K[X]$ (multiplication by $f_1(y)\in B^+$ acts on $V$), and $P_{k+1}-P_{k-1}=XP_k-2P_{k-1}$ by the recurrence of the $P_m$. Hence
\[
\Lambda\,\iota(f)=\tfrac12\vartheta\,\iota(f)=\iota\Bigl(\bigl(\tfrac X2P_k-P_{k-1}\bigr)f\Bigr)=\iota(-P_{k-1}f):
\]
under $\iota$ the operator $\Lambda$ on $V$ corresponds to multiplication by $-P_{k-1}$ on $R=K[X]/(P_k)$.
\item $V$ embeds into the quotient ring $B/QB$, since the kernel of the map $B^-\to B/QB$ is $QB\cap B^-$. On the image of $V$ in $B/QB$, $\Lambda$ agrees with multiplication by $t^k(t-x)$, and also with multiplication by $t^{-k}(t^{-1}-x)$: indeed $t^k(t-x)-t^{-k}(t^{-1}-x)=t^{-k-1}Q(t)\in QB$, and $\vartheta/2$ is the half-sum of these two elements. This description is valid only modulo $Q$: multiplication by $t^k(t-x)$ does not map $B^-$ into $B^-$.
\end{itemize}
\end{remark}

\begin{remark}[Symplectic characters]\label{rem:symplectic}
Let $\rho_1,\dots,\rho_k$ be the roots of $P_k$ in an algebraic closure $\Omega$ of $K$, and choose $\omega_j\in\Omega$ with $\omega_j+\omega_j^{-1}=\rho_j$. The $\rho_j$ are pairwise distinct (the discriminant of $P_k$ is a polynomial in $x$ which does not vanish at $x=0$, where $P_k=U_k$ has the distinct roots $2\cos\frac{j\pi}{k+1}$), and $\rho_j\ne\pm2$ because $P_k(2)=k+1-kx$ and $P_k(-2)=(-1)^k(k+1+kx)$ are non-zero; hence the $2k$ elements $\omega_j^{\pm1}$ are pairwise distinct. Let $C_m$ be the coordinate matrix in the definition of $w(m)$ (Theorem~\ref{thm:C}), so that $w(m)=\det C_m$. Since $P_k(\rho_j)=0$, evaluation at $\rho_1,\dots,\rho_k$ gives $\bigl(U_{m-1+i}(\rho_j)\bigr)_{i,j}=C_m\bigl(U_l(\rho_j)\bigr)_{l,j}$, and $\det\bigl(U_l(\rho_j)\bigr)_{0\le l\le k-1,\,1\le j\le k}=\prod_{i<j}(\rho_j-\rho_i)\ne0$ because $U_l$ is monic of degree $l$. Hence
\[
w(m)=\frac{\det\bigl(U_{m-1+i}(\rho_j)\bigr)_{0\le i\le k-1,\,1\le j\le k}}{\det\bigl(U_i(\rho_j)\bigr)_{0\le i\le k-1,\,1\le j\le k}}.
\]
By \eqref{eq:cheb} with $t\mapsto\omega_j$, $(\omega_j-\omega_j^{-1})U_n(\rho_j)=\omega_j^{n+1}-\omega_j^{-n-1}$, so by the Weyl character formula for $\Sp(2k)$ \cite[\S24.2, (24.18)]{FH}
\[
w(m)=\frac{\det\bigl(\omega_j^{m+i}-\omega_j^{-m-i}\bigr)}{\det\bigl(\omega_j^{1+i}-\omega_j^{-1-i}\bigr)}=\symp_{((m-1)^k)}(\omega_1,\dots,\omega_k)\qquad(m\ge1),
\]
the character of the irreducible representation of $\Sp(2k)$ indexed by the rectangular partition $((m-1)^k)$, evaluated at $\omega_1,\dots,\omega_k$. Moreover $x=\rho_1+\dots+\rho_k=\symp_{(1)}(\omega_1,\dots,\omega_k)$, because the coefficient of $X^{k-1}$ in $P_k=U_k-xU_{k-1}$ is $-x$. Thus for $\ell\ge k+2$ the identity \eqref{eq:C2} is a bilinear relation between rectangular symplectic characters at the point $(\omega_1,\dots,\omega_k)$, where $\omega_1^{\pm1},\dots,\omega_k^{\pm1}$ are the roots of $Q(t)/(t^2-1)=t^kP_k(t+t^{-1})$. It is not a relation at a generic point: if $x$ is replaced by $\symp_{(1)}$, the relation fails numerically at randomly chosen points $(\omega_1,\dots,\omega_k)$ already for $k=2,3$. We do not know whether \eqref{eq:C2} is a specialization of a known identity for symplectic characters. In \cite{Sun}, \eqref{eq:C2} is obtained from an identity with a free parameter between symplectic characters that holds at every point, by fixing the parameter and restricting to the points with $P_k=U_k-xU_{k-1}$.
\end{remark}

\subsection*{Computational checks}
The proofs in this paper do not depend on any computation. As an independent check, the statements and the main intermediate steps were verified by computer in exact arithmetic, using several independent implementations: Theorem~\ref{thm:A} with $x$ and $c$ indeterminates for $k\le10$ and $n\le22$; Theorem~\ref{thm:B}, with $D_{k,n}$ computed from its definition, in $\Z[x,c]$ for $k\le8$ and $2k+2\le n\le30$, with $x$ an indeterminate and $c\in\{0,3,-5\}$ for $k\le20$ and $2k+2\le n\le2k+42$, and in the form for $d_{k,n}(x)=\det H_{k,n}(x)$ for $k\le10$ and $2k+2\le n\le2k+30$; both formulas of Theorem~\ref{thm:C} with $x$ an indeterminate for $k\le14$, the first for $0\le n\le2k+32$ and the second for $-50\le\ell\le50$; and every step of the proof of Proposition~\ref{prop:bilinear} at several rational values of $x$ for $k\le6$ and all windows with $-3k-4\le\ell<3k+20$, including all degenerate windows in this range. The values of $D_{1,n}$, $D_{2,n}$ and $D_{3,n}$ computed from the definition agree with the lists in \cite{CMO} (and, for $D_{1,n}$ and $D_{2,n}$, with those on p.~5 of \cite{C23}), and the formulas of Remark~\ref{rem:smallk} were compared with the definition of $D_{k,n}$ for $n\le56$.

\section*{Acknowledgments}
The author thanks Haiming Sun for a careful reading of the manuscript and suggestions on its presentation.

\appendix

\section{Proof of Theorem~\ref{thm:A}}\label{app:A}

Let $J_0$ be the matrix $J$ of Section~\ref{sec:prelim} with $x$ replaced by $0$ (so that all diagonal entries of $J_0$ are equal to $c$), and let $a^\circ_{n,j}:=a_{n,j}(0,c)$; this is the triangle with $s_j=c$ for all $j$. Thus $J=J_0+xE_{00}$, where $E_{00}=e_0e_0^{\T}$, and by Lemma~\ref{lem:transfer} (specialized to $x=0$) $a^\circ_{n,j}=(J_0^n)_{0,j}$, $a^\circ_{n,j}=0$ for $j>n$ and $a^\circ_{n,n}=1$.

\begin{lemma}\label{lem:cheb-lin}
\begin{enumerate}
\item[(i)] $U_m(J_0-cI)\,e_0=e_m$ for all $m\ge0$.
\item[(ii)] For integers $p\ge k\ge0$, $U_pU_k=\sum_{i=0}^kU_{p-k+2i}$.
\end{enumerate}
\end{lemma}

\begin{proof}
(i) We have $(J_0-cI)e_0=e_1$ and $(J_0-cI)e_m=e_{m+1}+e_{m-1}$ for $m\ge1$. The claim holds for $m=0,1$, and inductively, for $m\ge1$, $U_{m+1}(J_0-cI)e_0=(J_0-cI)e_m-e_{m-1}=e_{m+1}$.

(ii) The ring homomorphism $\Z[X]\to\Z[t,t^{-1}]$, $X\mapsto t+t^{-1}$, is injective (the image of a polynomial of degree $d$ has highest exponent $d$), and $\Z[t,t^{-1}]$ is a domain. By \eqref{eq:cheb}, $U_k(t+t^{-1})=\sum_{i=0}^kt^{k-2i}$, since $(t-t^{-1})\sum_{i=0}^kt^{k-2i}=t^{k+1}-t^{-k-1}$. Hence, evaluating at $X=t+t^{-1}$,
\begin{align*}
(t-t^{-1})\,U_pU_k&=\sum_{i=0}^kt^{k-2i}\bigl(t^{p+1}-t^{-p-1}\bigr)=\sum_{i=0}^k\bigl(t^{p+k-2i+1}-t^{-(p+k-2i+1)}\bigr)\\
&=(t-t^{-1})\sum_{i=0}^kU_{p+k-2i}(t+t^{-1})
\end{align*}
(the second sum $\sum_it^{k-2i-p-1}$ is reindexed by $i\mapsto k-i$, and the last step is \eqref{eq:cheb} with $m=p+k-2i\ge0$). Cancelling $t-t^{-1}$ and using injectivity gives $U_pU_k=\sum_{i=0}^kU_{p+k-2i}=\sum_{i=0}^kU_{p-k+2i}$.
\end{proof}

\begin{lemma}\label{lem:gf}
Let $\Gamma=\Gamma(z)\in\Z[c][[z]]$ be the unique formal power series with $\Gamma=1+cz\Gamma+z^2\Gamma^2$. Then for every $q\ge0$
\[
\varphi^\circ_q(z):=\sum_{j\ge0}a^\circ_{j,q}z^j=z^q\Gamma^{q+1}.
\]
\end{lemma}

\begin{proof}
The series $\Gamma$ exists and is unique, since comparing coefficients of $z^n$ in $\Gamma=1+z(c\Gamma+z\Gamma^2)$ expresses the $n$-th coefficient of $\Gamma$ through coefficients of index $<n$. By the recurrence $a^\circ_{j+1,q}=a^\circ_{j,q-1}+c\,a^\circ_{j,q}+a^\circ_{j,q+1}$, $a^\circ_{0,q}=\delta_{q,0}$, $a^\circ_{j,-1}=0$, the family $(\varphi^\circ_q)_{q\ge0}$ satisfies
\begin{equation}\label{eq:phi-rec}
\varphi^\circ_q=\delta_{q,0}+z\bigl(\varphi^\circ_{q-1}+c\,\varphi^\circ_q+\varphi^\circ_{q+1}\bigr)\qquad(q\ge0,\text{ with }\varphi^\circ_{-1}:=0),
\end{equation}
and a family of power series satisfying \eqref{eq:phi-rec} is unique, because \eqref{eq:phi-rec} fixes the constant terms and expresses the coefficient of $z^{j+1}$ of each member through coefficients of $z^j$. The family $(z^q\Gamma^{q+1})_{q\ge0}$ also satisfies \eqref{eq:phi-rec}: for $q\ge1$, $z\bigl(z^{q-1}\Gamma^q+cz^q\Gamma^{q+1}+z^{q+1}\Gamma^{q+2}\bigr)=z^q\Gamma^q\bigl(1+cz\Gamma+z^2\Gamma^2\bigr)=z^q\Gamma^{q+1}$, and for $q=0$, $1+z(c\Gamma+z\Gamma^2)=\Gamma$. Hence $\varphi^\circ_q=z^q\Gamma^{q+1}$.
\end{proof}

\begin{lemma}[Paths between arbitrary heights]\label{lem:arb}
For all $p,k,j\ge0$,
\[
(J^j)_{p,k}=\sum_{q\ge0}h^{(k)}_{p,q}\,a^\circ_{j,q}
\]
$($a finite sum, since $a^\circ_{j,q}=0$ for $q>j)$. Here, equivalently to the definition in \S1.1, $h^{(k)}_{p,q}=1$ if $q\in\{|p-k|,|p-k|+2,\dots,p+k\}$, $h^{(k)}_{p,q}=x^{q-p-k}$ if $q>p+k$, and $h^{(k)}_{p,q}=0$ otherwise.
\end{lemma}

\begin{proof}
The two descriptions of $h^{(k)}_{p,q}$ agree: $|p-q|\le k\le p+q$ means $|p-k|\le q\le p+k$, and $p+q\equiv k\pmod2$ means $q\equiv|p-k|\pmod2$. Put $\varphi_{p,k}(z):=\sum_{j\ge0}(J^j)_{p,k}z^j$ and $\varphi^\circ_{p,k}(z):=\sum_{j\ge0}(J_0^j)_{p,k}z^j$ in $\Z[x,c][[z]]$; thus $\varphi^\circ_{0,q}=\varphi^\circ_q$.

\emph{Step 1: $\varphi^\circ_{p,k}=\sum_{i=0}^{\min(p,k)}\varphi^\circ_{|p-k|+2i}$.} Both sides are symmetric in $(p,k)$ (because $J_0$ is symmetric), so let $p\ge k$. By Lemma~\ref{lem:cheb-lin}, the symmetry of $U_p(J_0-cI)$ and the commutativity of polynomials in $J_0$,
\begin{align*}
(J_0^j)_{p,k}&=\bigl(U_p(J_0-cI)e_0\bigr)^{\T}J_0^jU_k(J_0-cI)e_0=e_0^{\T}J_0^j(U_pU_k)(J_0-cI)e_0\\
&=\sum_{i=0}^ke_0^{\T}J_0^je_{p-k+2i}=\sum_{i=0}^ka^\circ_{j,p-k+2i}.
\end{align*}
In particular $\varphi^\circ_{p,0}=\varphi^\circ_{0,p}=\varphi^\circ_p$.

\emph{Step 2: $\varphi_{p,k}=\varphi^\circ_{p,k}+\sum_{i\ge1}x^i\varphi^\circ_{p+k+i}$.} In $M(\Z[x,c])$ we have $J^j-J_0^j=\sum_{i=0}^{j-1}J_0^i(J-J_0)J^{j-1-i}$ and $J-J_0=xe_0e_0^{\T}$, so $(J^j)_{p,k}=(J_0^j)_{p,k}+x\sum_{i=0}^{j-1}(J_0^i)_{p,0}(J^{j-1-i})_{0,k}$, that is, $\varphi_{p,k}=\varphi^\circ_{p,k}+xz\,\varphi^\circ_p\,\varphi_{0,k}$. For $p=0$ this reads $\varphi_{0,k}=\varphi^\circ_k+xz\Gamma\varphi_{0,k}$, so $\varphi_{0,k}=\varphi^\circ_k/(1-xz\Gamma)$. Substituting back and using Lemma~\ref{lem:gf},
\[
\varphi_{p,k}=\varphi^\circ_{p,k}+\frac{xz\,\varphi^\circ_p\varphi^\circ_k}{1-xz\Gamma}=\varphi^\circ_{p,k}+\sum_{i\ge1}x^iz^i\Gamma^{i-1}\cdot z^p\Gamma^{p+1}\cdot z^k\Gamma^{k+1}=\varphi^\circ_{p,k}+\sum_{i\ge1}x^i\varphi^\circ_{p+k+i}
\]
(the infinite sums converge in the $z$-adic topology, because $\varphi^\circ_q\in z^q\Z[c][[z]]$).

Steps 1 and 2 give $\varphi_{p,k}=\sum_{q\ge0}h^{(k)}_{p,q}\varphi^\circ_q$; taking the coefficient of $z^j$ proves the lemma.
\end{proof}

Combinatorially, Step 1 is the reflection principle (or Clebsch--Gordan rule) for paths from height $p$ to height $k$ without level steps of weight $x$, and the term $x^i\varphi^\circ_{p+k+i}$ counts the paths with exactly $i$ level steps of weight $x$ at height $0$.

\begin{proof}[Proof of Theorem~\ref{thm:A}]
Let $n\ge1$ and $0\le i,j\le n-1$. By Lemma~\ref{lem:transfer}, associativity in $M(\Z[x,c])$, and Lemma~\ref{lem:arb},
\[
a_{i+j,k}=(J^{i+j})_{0,k}=\sum_{p\ge0}(J^i)_{0,p}(J^j)_{p,k}=\sum_{p\ge0}\sum_{q\ge0}a_{i,p}\,h_{p,q}\,a^\circ_{j,q}.
\]
Since $a_{i,p}=0$ for $p>i$ and $a^\circ_{j,q}=0$ for $q>j$, only $p,q\le n-1$ contribute. Hence
\[
\bigl(a_{i+j,k}\bigr)_{0\le i,j\le n-1}=\mathcal A_n\cdot H_{k,n}(x)\cdot(\mathcal A^\circ_n)^{\T},\qquad \mathcal A_n:=\bigl(a_{i,p}\bigr)_{0\le i,p\le n-1},\quad \mathcal A^\circ_n:=\bigl(a^\circ_{j,q}\bigr)_{0\le j,q\le n-1}.
\]
By Lemma~\ref{lem:transfer}, $\mathcal A_n$ and $\mathcal A^\circ_n$ are lower triangular with unit diagonal, so $\det\mathcal A_n=\det\mathcal A^\circ_n=1$ and $D_{k,n}=\det H_{k,n}(x)$. For $n=0$ both sides are equal to $1$.
\end{proof}

\end{document}